\documentclass[12pt]{article}
\usepackage{amsmath,amssymb, amsthm, graphicx}
\usepackage{color}
\newtheorem{thm}{Theorem}
\newtheorem{conj}[thm]{Conjecture}

\newcommand{\partit}{{\cal Z}}
\newcommand{\config}{{\cal C}}

\title{A formula for the doubly refined enumeration of alternating sign matrices}
\author{Matan Karklinsky \and Dan Romik
\thanks{Supported by the Israel Science Foundation
(ISF) grant number 1051/08.}
}
\begin{document}
\maketitle

\abstract{Zeilberger \cite{zeilberger2} proved the Refined Alternating Sign Matrix Theorem, which gives a product formula, first conjectured by Mills, Robbins and Rumsey \cite{millsetal2}, for the number of alternating sign matrices with given top row. Stroganov \cite{stroganov} proved an explicit formula for the number of alternating sign matrices with given top and bottom rows. Fischer and Romik \cite{fischerromik} considered a different kind of ``doubly-refined enumeration'' where one counts alternating sign matrices with given top two rows, and obtained partial results on this enumeration. In this paper we continue the study of the doubly-refined enumeration with respect to the top two rows, and  use Stroganov's formula to prove an explicit formula for  these doubly-refined enumeration numbers.
}

\section{Introduction}

An \textbf{alternating sign matrix (ASM)} of order $n$ is an $n\times n$ matrix with entries in $\{0,-1,1\}$ such that in every row and every column, the sum of the entries is $1$ and the non-zero numbers appear with alternating signs. See Fig.~\ref{fig:asmexample}(a) for an example. A \textbf{monotone triangle} of order $n$
is a triangular array
$ (t_{i,j})_{1\le i\le n, 1\le j\le i} $ of integers satisfying the inequalities
$$ t_{i,j} < t_{i,j+1}, \ \ \  t_{i,j} \le t_{i-1,j} \le t_{i,j+1} \qquad (2\le i\le n, 1\le j\le i-1). $$
A monotone triangle of order $n$ is called \textbf{complete} if its bottom row has the numbers $(1,2,\ldots,n)$. See 
Fig.~\ref{fig:asmexample}(b).

For integers $k_1<k_2<\ldots<k_n$, denote by $\alpha_n(k_1,k_2,\ldots, k_n)$ the number of monotone triangles with bottom row $(k_1,\ldots,k_n)$. Define
\begin{eqnarray*}
A_n &=& \alpha_n(1,2,\ldots,n), \\
A_{n,k} &=& \alpha_{n-1}(1,2,\ldots,\hat{k},\ldots,n), \ \,\quad\qquad (1\le k\le n)  \\
A_{n,i,j} &=& \alpha_{n-2}(1,\ldots,\hat{i},\ldots,\hat{j},\ldots,n), \qquad (1\le i<j\le n).
\end{eqnarray*}
where the notation $\hat{s}$ in a list of numbers indicates that $s$ is omitted from the list. For notational convenience later on, take $A_{n,k}$ as $0$ if $k\notin[1,n]$. 

Since ASMs of order $n$ are well-known to be in bijection with complete monotone triangles of order $n$, the number $A_n$ is the total number of ASMs of order $n$. From obvious properties of the bijection, it follows that $A_{n,k}$ is the number of ASMs of order $n$ such that the unique $1$ in the first row is in position $k$. The numbers $(A_{n,k})_{n,k}$ are called the refined enumeration numbers for ASMs. Zeilberger \cite{zeilberger1,zeilberger2} proved the celebrated alternating sign matrix theorem and refined alternating sign matrix theorem, which state, respectively, that
\begin{eqnarray}
A_n &=& \prod_{j=0}^{n-1} \frac{(3j+1)!}{(n+j)!}\qquad\text{and} \label{eq:asmthm} \\
A_{n,k} &=& \binom{n+k-2}{k-1} \frac{(2n-k-1)!}{(n-k)!}\prod_{j=0}^{n-2} \frac{(3j+1)!}{(n+j)!}. \label{eq:refinedasmthm}
\end{eqnarray}

\begin{figure}
\begin{center}
\begin{tabular}{cc} 
$  \left( \begin{array}{cccccc}
0 & 0 & 1 & 0 & 0 & 0 \\
0 & 1 & -1 & 0 & 1 & 0 \\
1 & -1 & 0 & 1 & 0 & 0 \\
0 & 1 & 0 & 0 & -1 & 1  \\
0 & 0 & 1 & 0 & 0 & 0 \\
0 & 0 & 0 & 0 & 1 & 0 
\end{array} \right)  
$
&
\begin{tabular}{*{15}{c@{\hspace{0.06in}}}}
& & && & 3 & & & & &
 \\
& & && 2 & & 5 & & & &
 \\
& && 1 & & 4 & & 5 & & &
 \\
&  & 1 && 2 && 4 & & 6 & &
 \\
& 1 && 2 && 3 && 4 && 6 &
\\
1 && 2 && 3 && 4 && 5 && 6
\end{tabular}
\\ 
(a)&(b) 
\end{tabular}
\caption{An ASM of order 5 and the corresponding complete monotone triangle.
\label{fig:asmexample}}
\end{center}
\end{figure}

In \cite{fischerromik} the numbers $(A_{n,i,j})_{n,i,j}$ were studied, in an attempt to extend the work of Zeilberger to a \textbf{doubly-refined enumeration} of alternating sign matrices based on their first two rows. We refer to these numbers as the doubly-refined enumeration numbers. The connection to enumeration of alternating sign matrices (discussed at length in \cite{fischerromik}) is illustrated in 
Fig.~\ref{fig:asmfirsttworows}, which shows the possible configurations for the first two rows of an ASM. For fixed $n$, these configurations are indexed by a triple $(i,j,k)$ such that $1\le i\le k\le j \le n$ and $i<j$, and for each such triple the number of ASMs having the given first two rows is independent of $k$ and in fact is equal to $A_{n,i,j}$.

\begin{figure}
\begin{center}
\begin{tabular}{ccc}
$ 
\begin{array}{l}
\begin{array}{cccccc}
\qquad \ i
\qquad \ \ k
\quad \ j
\end{array}\\ 
\left(\begin{array}{cccccc}
0 & 0 & 0 & 1 & 0 & 0 \\
0 & 1 & 0 & -1 & 1 & 0 \\
& & & \vdots & & 
\end{array}
\right)
\end{array}
$
&
$ 
\begin{array}{l}
\begin{array}{cccccc}
\qquad \qquad\ i\!\!=\!\!k
\  \ j
\end{array}\\ 
\left(\begin{array}{cccccc}
0 & 0 & 0 & 1 & 0 & 0 \\
0 & 0 & 0 & 0 & 1 & 0 \\
& & & \vdots & & 
\end{array}
\right)
\end{array}
$
&
$ 
\begin{array}{l}
\begin{array}{cccccc}
\qquad \ i
\quad\ \  k\!\!=\!\!j
\end{array}\\ 
\left(\begin{array}{cccccc}
0 & 0 & 0 & 1 & 0 & 0 \\
0 & 1 & 0 & 0 & 0 & 0 \\
& & & \vdots & & 
\end{array}
\right)
\end{array}
$
\\ (a) & (b) & (c)
\end{tabular}
\caption{The possible configurations for the first two rows of an ASM: (a) $i<k<j$. (b) $i=k<j$. (c) $i<k=j$. In all three cases the number of ASMs with the given first two rows is 
$A_{n,i,j}$.\label{fig:asmfirsttworows}}
\end{center}
\end{figure}

The results of \cite{fischerromik} gave only partial information on the $A_{n,i,j}$'s, namely a system of linear equations satisfied by
$(A_{n,i,j})_{i,j}$ for each $n$. These equations, along with some other information known from simple considerations, were conjectured to determine the $A_{n,i,j}$'s uniquely and so to allow to express the $A_{n,i,j}$'s as ratios of determinants and to compute them efficiently on a computer for reasonably large values of $n$. A rather complicated explicit formula not involving determinants was also conjectured based on numerical evidence.

Our main result is a new and simpler explicit formula for $A_{n,i,j}$. We prove:

\begin{thm} Let
\begin{equation}\label{eq:def-xnij}
X_n(s,t) = \frac{1}{A_{n-1}} \Bigg(A_{n-1,t}(A_{n,s+1}-A_{n,s})-A_{n-1,s}(A_{n,t+1}-A_{n,t}) \Bigg).
\end{equation}
For each  $1\le i<j\le n$ we have
\begin{equation}\label{eq:explicit-formula}
A_{n,i,j} =  \sum_{p=0}^{n-j} \sum_{q=0}^p (-1)^q \binom{p}{q} X_n(i+q,j+p).
\end{equation}
\label{thm:main}
\end{thm}

Our proof of Theorem \ref{thm:main} uses the well-known connection between alternating sign matrices and the \textbf{square ice} model, and builds on previous results and techniques of Izergin-Korepin \cite{bogoliubov-izergin-korepin}, Kuperberg \cite{kuperberg} and Stroganov \cite{stroganov}. Following our discovery of \eqref{eq:explicit-formula}, an alternative derivation of the same formula using the monotone triangle techniques developed in \cite{fischer1, fischer2, fischerromik} was recently found by Fischer~\cite{fischerprivate}.

The following conjecture holds empirically but does not follow from our methods.

\begin{conj}
For each $n\ge 1$ let $\left(\hat{A}_{n,i,j}\right)_{i,j=1}^n$ be the \emph{extended} doubly-refined enumeration numbers defined in \cite{fischerromik}. Then \eqref{eq:explicit-formula} holds also for the extended numbers. In other words, we have
$$
\hat{A}_{n,i,j} =  \sum_{p=0}^{n-j} \sum_{q=0}^p (-1)^q \binom{p}{q} X_n(i+q,j+p), \qquad (1 \le i,j \le n).
$$
\end{conj}

\section{Square ice and the partition function}

\begin{figure}
\begin{center}
\begin{tabular}{cccccc}
\scalebox{0.2}{\includegraphics{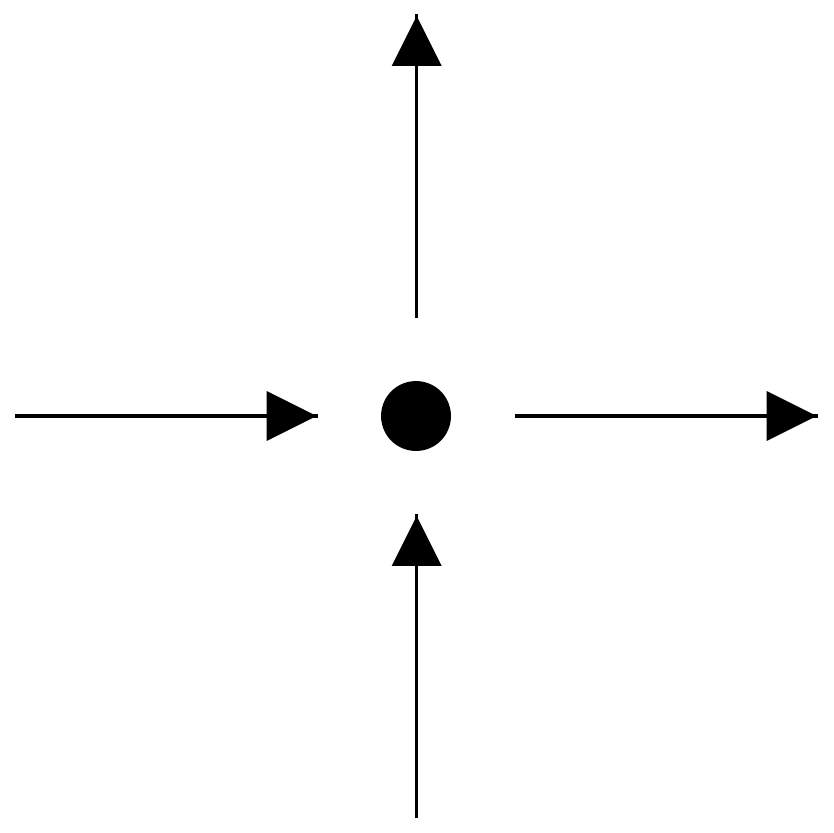}} &
\scalebox{0.2}{\includegraphics{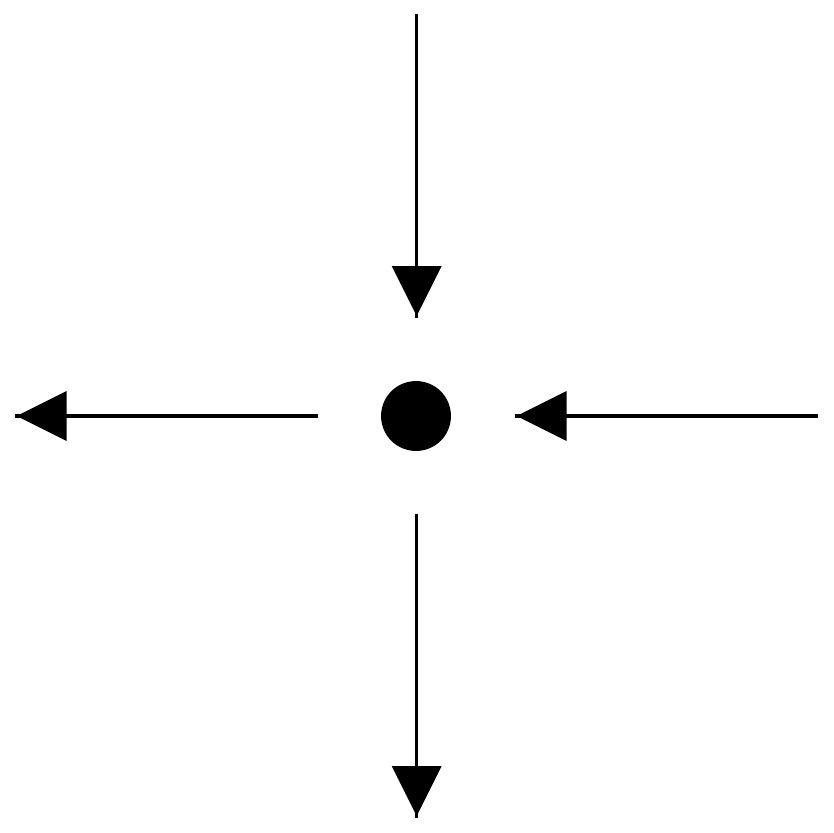}} &
\scalebox{0.2}{\includegraphics{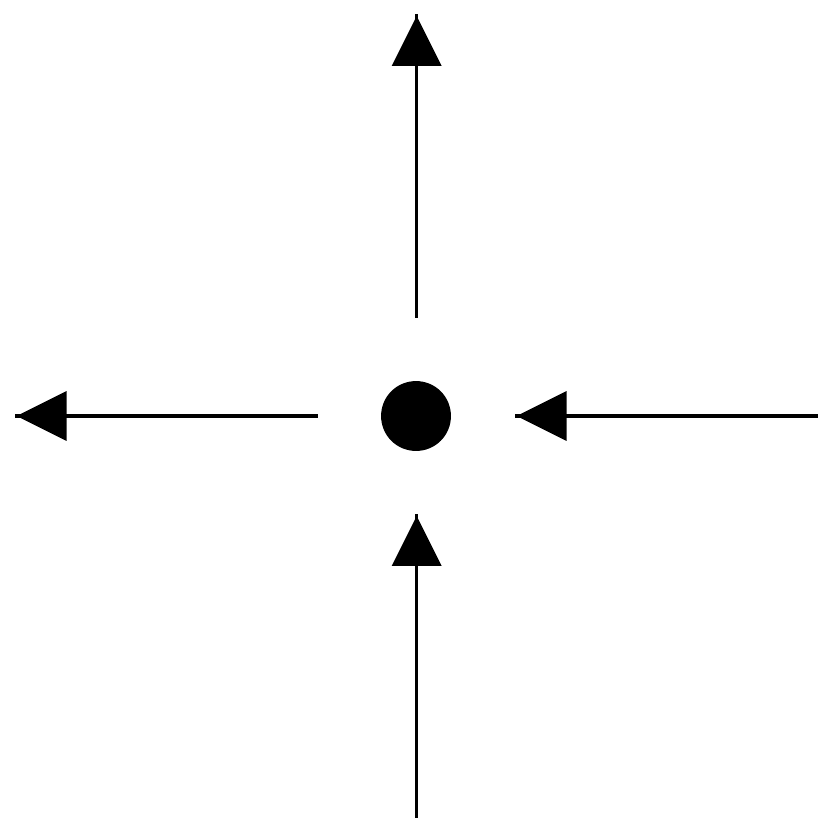}} &
\scalebox{0.2}{\includegraphics{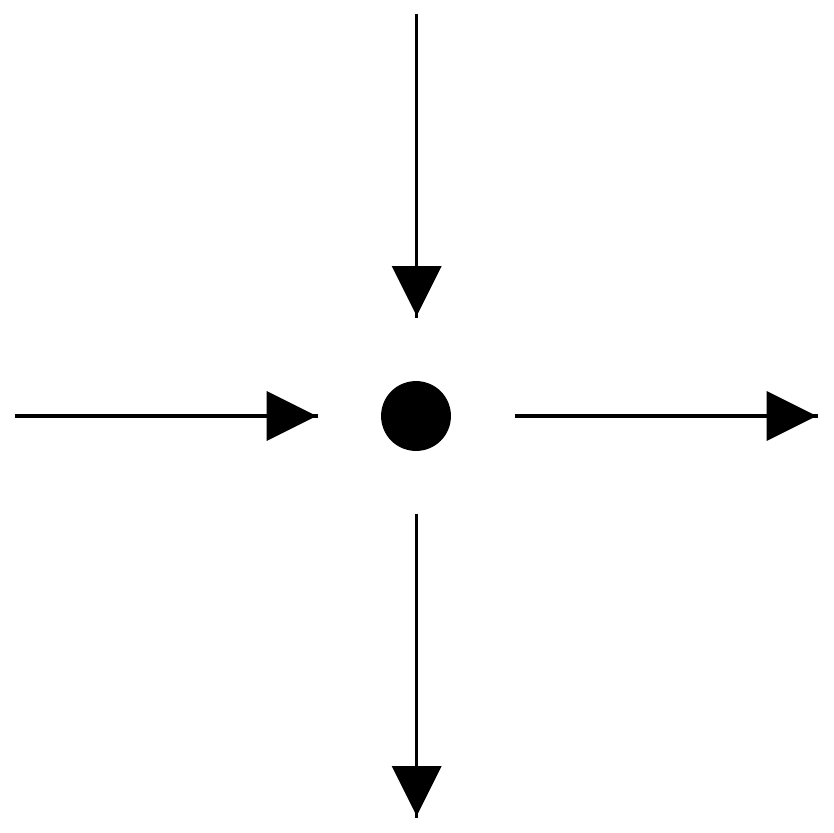}} &
\scalebox{0.2}{\includegraphics{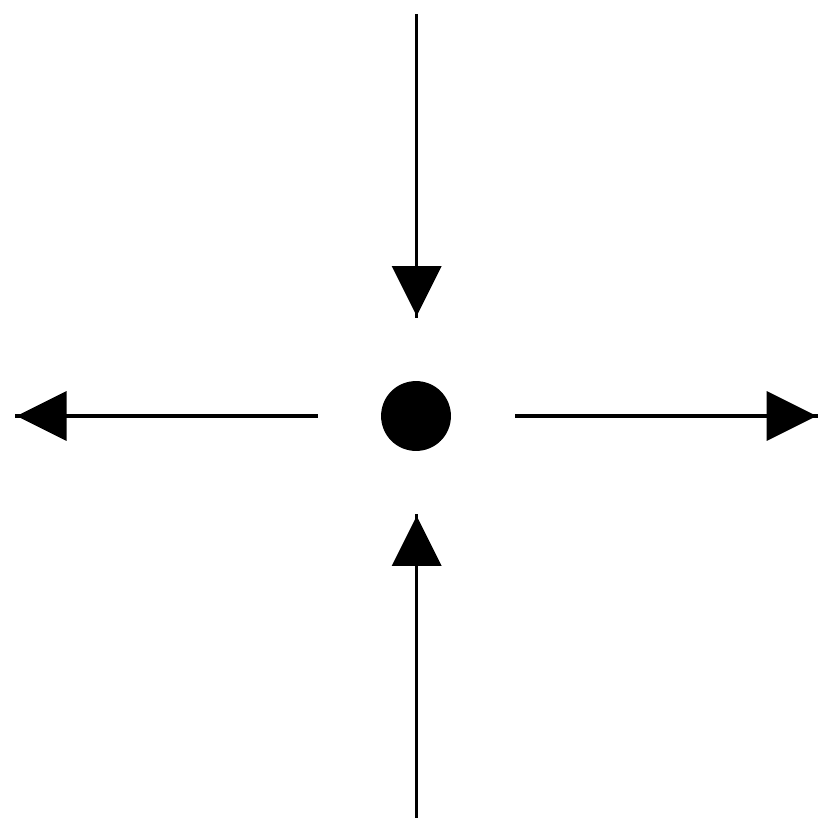}} &
\scalebox{0.2}{\includegraphics{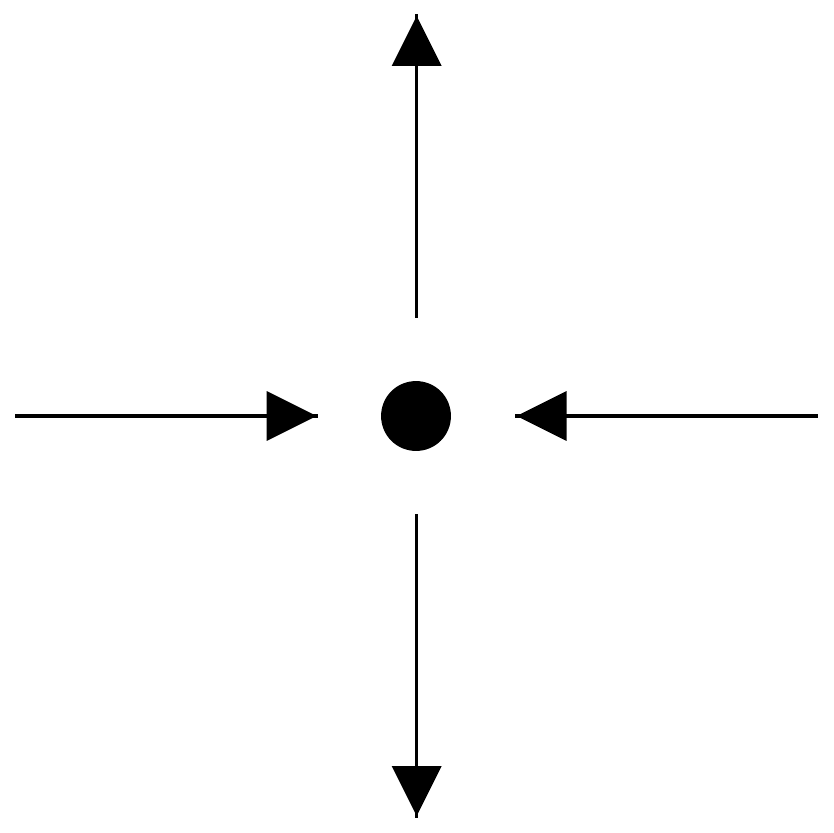}} \\
(type 1) &
(type 2) &
(type 3) &
(type 4) &
(type 5) &
(type 6)
\\
$w(v)=a$ &
$w(v)=a$ &
$w(v)=b$ &
$w(v)=b$ &
$w(v)=c$ &
$w(v)=c$ \\
$m = 0$ &
$m = 0$ &
$m = 0$ &
$m = 0$ &
$m = -1$ &
$m = 1$
\end{tabular}
\caption{The six types of vertices in a square ice configuration, the weight $w$ attached to each of them and the entry $m$ assigned in the corresponding ASM under the bijection translating square ice configurations to ASM's.
\label{fig:sixvertex}}
\end{center}
\end{figure}

\begin{figure}
\begin{center}
\scalebox{0.75}{\includegraphics{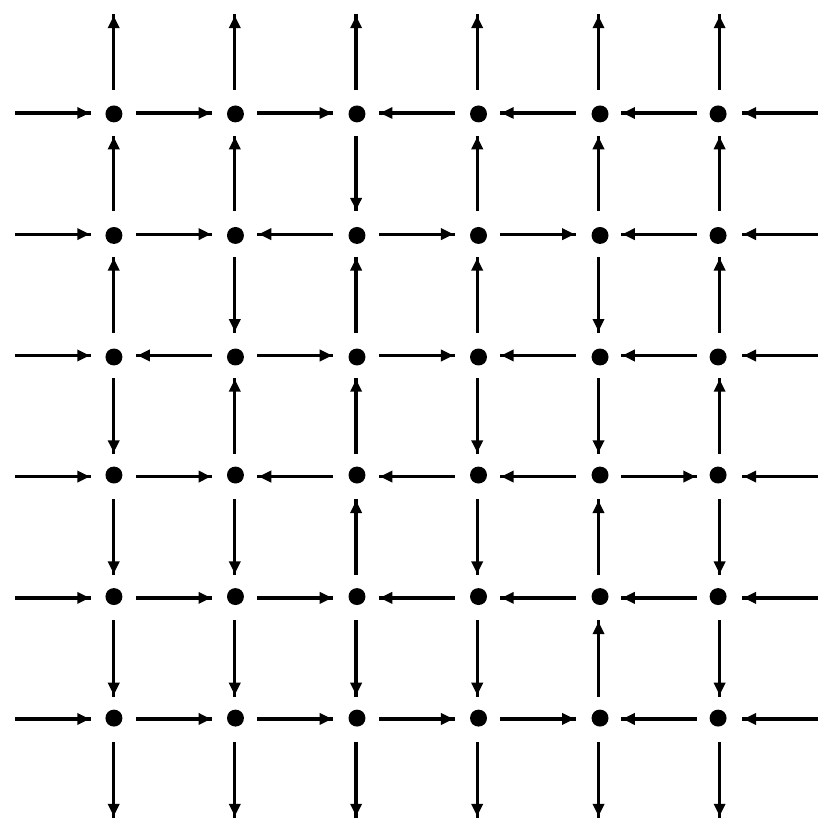}}
\caption{The square ice configuration corresponding to the ASM in Figure \ref{fig:asmexample}(a)\label{fig:asmice-bijection}.}
\end{center}
\end{figure}

We consider square ice (or \textbf{six-vertex model}) configurations on an $n\times n$ square lattice 
$$L_n=\{1,\ldots,n\}\times\{1,\ldots,n\}$$ 
satisfying the so-called \textbf{domain wall boundary conditions}. In a square ice configuration the edges of the lattice are oriented so that each vertex in the lattice has two incoming edges and two outgoing edges, giving six possibilities, shown in Fig.~\ref{fig:sixvertex}. The boundary conditions are that edges entering the lattice from the left and right sides of the square point inwards, whereas edges adjacent to the top and bottom sides point outwards. Such configurations are in bijection with alternating sign matrices of order $n$. The bijection maps each of the six types of square ice vertices to either a $+1$, a $-1$ or a $0$; see Figs. \ref{fig:sixvertex} and \ref{fig:asmice-bijection}.

Fix  a real value $\eta$ called the \textbf{crossing parameter}. Given a square ice configuration $C$, and given a vertex $v$ in the lattice with an associated (real-valued) \textbf{row parameter} $x$ and \textbf{column parameter} $y$, according to the type of the vertex $v$ in the configuration $C$ we associate with it a \textbf{weight} $w(v)=w_C(v)$ equal to either of the three quantities $a, b$ or $c$, defined by
\begin{eqnarray*}
a &=& \frac{\sin(\eta/2+x-y)}{\sin\eta}, \\
b &=& \frac{\sin(\eta/2-x+y)}{\sin\eta}, \\
c &=& 1.
\end{eqnarray*}
Vertices of type 1 and 2 are assigned the weight $a$, those of type 3 and 4 are assigned the weight $b$, and types 5 and 6 get the weight $c$.
The weight of a square ice configuration $C$ is then defined as the product of the vertex weights over all the vertices in the lattice, namely
$$ w(C) = \prod_{v\in L_n} w_C(v). $$
Denote the set of square ice configurations by $\config_n$.
The square ice \textbf{partition function} is defined as the sum of the weights over all configurations and is denoted by $\partit_n$, i.e.,
$$ \partit_n = \sum_{C\in \config_n} w(C) = \sum_{C\in \config_n} \prod_{v\in L_n} w_C(v). $$
Here, we associate with a vertex in row $i$ and column $j$ a row parameter $x_i$ and a column parameter $y_j$. Thus the partition function $\partit_n$ is a function of the crossing parameter $\eta$ and of the variables $x_1,\ldots,x_n,y_1,\ldots,y_n$, called the \textbf{spectral parameters}. Regarding $\eta$ as fixed, we occasionally write $\partit_n(x_1,\ldots,x_n; y_1,\ldots,y_n)$ to emphasize the dependence on the spectral parameters or to substitute specific values.

An important fact that we will use is that $\partit_n$ is a symmetric function in the row parameters $x_1,\ldots, x_n$ and is also symmetric in the column parameters $y_1,\ldots,y_n$. For the remainder of the paper we take $\eta=2\pi/3$. For this particular value of the crossing parameter, Stroganov \cite{stroganov} showed that $\partit_n$ is actually symmetric in the union of the $x_i$'s and $y_j$'s, but we will not need this fact.

For our proof of Theorem \ref{thm:main}, we will consider two specializations of the partition function, namely
\begin{eqnarray*}
S_1 &=& \partit_n(u,0,\ldots,0,v ; 0,\ldots,0) \textrm{ and}\\
S_2 &=& \partit_n(u,v,0,\ldots,0 ; 0,\ldots,0)
\end{eqnarray*}
where $u$ and $v$ are parameters.
By the symmetry of $\partit_n$ in the row parameters, we get that $S_1=S_2$. This will enable us to relate the doubly-refined enumeration numbers to a different kind of doubly-refined enumeration studied by Stroganov, where one enumerates ASM based on the first and last row instead of the first two rows.

\section{Evaluation of $S_1$}

For $n\ge 1$ and $1\le i,j\le n$, let $B_{n,i,j}$ denote the number of ASM's of order $n$ with a $1$ in positions $(1,i)$ and $(n,j)$. These numbers enumerate ASM's based on their first and last row and so are also sometimes referred to as the doubly-refined enumeration numbers. They were originally considered by Mills, Robbins and Rumsey \cite{millsetal2}. Stroganov proved the following formula expressing $B_{n,i,j}$ in terms of the (singly-) refined enumeration numbers $A_{n,k}$
(see also \cite{colomopronkorefinedasm} where some generalizations of this result are proved).

\begin{thm}[Stroganov \cite{stroganov}] We have 
$$
B_{n,i+1,j+1} - B_{n,i,j} =  Y_n(i,j) := 
\frac{1}{A_{n-1}} \Bigg(A_{n-1,j}(A_{n,i+1}-A_{n,i})+A_{n-1,i}(A_{n,j+1}-A_{n,j}) \Bigg).
$$
\end{thm}

Note that $Y_n(i,j) = X_n(i,n-j)$ where $X_n(s,t)$ is defined in \eqref{eq:def-xnij}. By summation of the differences it follows also that
$$
B_{n,i,j} = A_{n-1,|i-j|} +  \sum_{k=1}^{\min(i,j)-1} Y_n(k,|i-j|+k).
$$
Stroganov's proof involved evaluating the expression denoted above by $S_1$ in two ways. For our purposes, we only need one of them, whose derivation we include for completeness. Let 
$$ \varphi(x) = \frac{2}{\sqrt{3}}\sin(\pi/3+x), \quad t=t(u)=\frac{\varphi(u)}{\varphi(-u)}, \quad s=s(v)=\frac{\varphi(-v)}{\varphi(v)}. $$
Now note that when summing over all square ice configurations to evaluate $S_1=\partit_n(u,0,\ldots,0,v ; 0,\ldots,0)$, the weight of each configuration depends only on the positions of the $1$'s in the first and last rows of the ASM corresponding to the configuration, since
all vertices in rows $2$ through $n-1$ contribute a factor of $1$.
If the ASM has 1's in positions $i$ and $j$ respectively in the first and last row, using the translation between ASM's and square ice configurations the weight of the corresponding configuration is easily seen to be
$$ \varphi(u)^{i-1} \varphi(-u)^{n-i} \varphi(-v)^{j-1} \varphi(v)^{n-j}. $$
This implies that
\begin{eqnarray}
S_1 &=&  \sum_{i,j=1}^n B_{n,i,j} \varphi(u)^{i-1} \varphi(-u)^{n-i} \varphi(-v)^{j-1} \varphi(v)^{n-j} \nonumber \\
& = &
\Big(\varphi(v)\varphi(-u)\Big)^{n-1} \sum_{i,j=1}^n B_{n,i,j} t^{i-1} s^{j-1}.
\label{eq:s1-eval}
\end{eqnarray}

\section{Evaluation of $S_2$}

To evaluate $S_2$, note again that for this substitution of spectral parameters the weight of a configuration is only dependent on the state of the first two rows of the square ice configuration, or equivalently of the corresponding ASM. As explained above, this state can be indexed by parameters $i,j,k$ satisfying $1\le i\le k\le j\le n$ and $i<j$, which correspond to the three numbers in the top two rows of the complete monotone triangle that corresponds to the ASM. The weight of the corresponding square ice configuration having given parameters $i,j,k$ is then easily computed to be
$$ w(C) = \begin{cases} \varphi(u)^{k-1}\varphi(-u)^{n-k} \varphi(v)^{i+j-k-2} \varphi(-v)^{n-j+k-i-1} & i<k<j, \\
\varphi(u)^{k-1} \varphi(-u)^{n-k} \varphi(v)^{j-2} \varphi(-v)^{n-j+1} & i=k<j, \\
\varphi(u)^{k-1} \varphi(-u)^{n-k} \varphi(v)^i \varphi(-v)^{n-i-1} & i<k=j.
 \end{cases}
 $$
Now sum this over all configurations, and use the ``120-degree triangle'' identity
$$ \varphi(x)^2 + \varphi(-x)^2 - \varphi(x)\varphi(y) = 1, $$
to get that
\begin{eqnarray*}
S_2 &=& \sum_{1\le i<j\le n} A_{n,i,j}
\Bigg(
\varphi(u)^{k-1} \varphi(-u)^{n-k} \varphi(v)^{j-2} \varphi(-v)^{n-j+1}
\\ & & \qquad\qquad\qquad\ \  + \varphi(u)^{k-1} \varphi(-u)^{n-k} \varphi(v)^i \varphi(-v)^{n-i-1}
\\ & & \qquad\qquad\qquad\ \  + \sum_{k=i+1}^{j-1}
\varphi(u)^{k-1}\varphi(-u)^{n-k} \varphi(v)^{i+j-k-2} \varphi(-v)^{n-j+k-i-1}
\Bigg)
\\ & &
= \Big(\varphi(v)\varphi(-u)\Big)^{n-1}
\sum_{1\le i<j\le n} A_{n,i,j}
\Bigg(
t^{i-1} s^{n-j+1}+t^{j-1}{n-i-1}
\\ & & \qquad\qquad\qquad\qquad\qquad\qquad\qquad\quad\  + \frac{1}{\varphi(v)^2} \sum_{k=i+1}^{j-1} t^{k-1} s^{n-j+k-i-1}
\Bigg)
\\ & &
= \Big(\varphi(v)\varphi(-u)\Big)^{n-1}
\sum_{1\le i<j\le n} A_{n,i,j}
\Bigg(
t^{i-1} s^{n-j+1}+t^{j-1}{n-i-1}
\\ & & \qquad\qquad\qquad\qquad\qquad\qquad\qquad\quad\  + (1+s^2-s) \sum_{k=i+1}^{j-1} t^{k-1} s^{n-j+k-i-1}
\Bigg).
\end{eqnarray*}
The inner summation on $k$ is a finite geometric series that sums to
$$
\sum_{k=i+1}^{j-1} t^{k-1} s^{n-j+k-i-1} = \frac{t^{j-1} s^{n-i-1} - t^i s^{n-j}}{t s - 1}, $$
so that, after some further simple algebraic simplifications, we obtain
\begin{eqnarray}
\frac{(ts - 1)S_2}{(\varphi(v)\varphi(-u))^{n-1}} &=&
\sum_{1\le i<j \le n}
A_{n,i,j} \Big(
t^i s^{n-j+1} -t^i s^{n-j} - t^{i-1} s^{n-j+1} 
\nonumber \\ & & 
\hspace{70.0pt}
+ \, t^j s^{n-i} + t^{j-1} s^{n-i+1} - t^{j-1} s^{n-i}
\Big).
\label{eq:s2-eval}
\end{eqnarray}

\section{Completion of the proof}

Now knowing that  $S_1 = S_2$ we can write
$$
\frac{(ts - 1)S_1}{(\varphi(v)\varphi(-u))^{n-1}} = \frac{(ts - 1)S_2}{(\varphi(v)\varphi(-u))^{n-1}}.
$$
We can equate the coefficients of $t^i s^j$ on both sides (since these functions are linearly independent, a fact that is equivalent to \cite[Ex. 7.1.9, p. 231]{bressoud}). In conjunction with \eqref{eq:s1-eval} and \eqref{eq:s2-eval} this translates to the identity
\begin{multline}
Y_n(i,j) = B_{n,i+1,j+1}-B_{n,i,j}
\\ = A_{n,i+1,n+1-j} + A_{n,i,n-j} - A_{n,i,n+1-j} - A_{n,n-j,i}
- A_{n,n+1-j, i+1} + A_{n,n-j,i+1},
\label{eq:ab-identity}
\end{multline}
which holds for all $i,j$ if we adopt the convention that $A_{n,p,q}$ and $B_{n,p,q}$ are $0$ outside the respective ranges of their definitions.

It remains to solve this system of equations in the $A_{n,i,j}$'s. First we reformulate it slightly for convenience by replacing $j$ by $n-j$. This gives
\begin{eqnarray*}
X_n(i,j) &=& Y_n(i,n-j)=
B_{n,i+1,n-j+1}-B_{n,i,n-j} \\ &=&
A_{n,i+1,j+1}+A_{n,i,j}-A_{n,i,j+1} - A_{n,j,i} - A_{n,j+1,i+1} + A_{n,j,i+1}.
\end{eqnarray*}
This simplifies even further when one observes that the last 3 terms on the right-hand side are 0 in the range of parameters $1\le i<j \le n$ which interests us. So we have
$$ X_n(i,j) = A_{n,i+1,j+1} + A_{n,i,j} - A_{n,i,j+1}, \qquad (1\le i<j\le n). $$
Rewriting this in the form
$$ A_{n,i,j} = X_n(i,j) + A_{n,i,j+1} - A_{n,i+1,j+1}, $$
it can now be easily solved by iteration (or more formally by reverse induction on $j$), as follows:
\begin{eqnarray*}
A_{n,i,j} &=& X_n(i,j)+A_{n,i,j+1}-A_{n,i+1,j+1}
\\ &=&
X_n(i,j)+X_n(i,j+1)-X_n(i+1,j+1)
\\ & &+(A_{n,i,j+2}-A_{n,i+1,j+2})-(A_{n,i+1,j+2}-A_{n,i+2,j+2})
\\ &=&
X_n(i,j)+X_n(i,j+1)-X_n(i+1,j+1) 
+A_{n,i,j+2}-2A_{n,i+1,j+2}+A_{n,i+2,j+2}
\\ &=&
X_n(i,j)+\big(X_n(i,j+1)-X_n(i+1,j+1)\big)
\\ & & +\big(X_n(i,j+2)-2 X_n(i+1,j+2)+X_n(i+2,j+2)\big)\\ & &
+A_{n,i,j+3}-3 A_{n,i+1,j+3} + 3 A_{n,i+2,j+3} - A_{n,i+3,j+3}
\\ &=& \ldots
\\ &=& \sum_{k=0}^{n-j} (-1)^k \binom{n-j}{k} A_{n,i+k,n} + \sum_{p=0}^{n-j-1} \sum_{q=0}^p \binom{p}{q} (-1)^q X_n(i+q,j+p).
\end{eqnarray*}
In the last expression, the first summation over $k$ corresponds exactly to the case $p=n-j$ of the second sum, so we can shorten this to
$$ A_{n,i,j} = \sum_{p=0}^{n-j} \sum_{q=0}^p \binom{p}{q} (-1)^q X_n(i+q,j+p) $$
which was the claim of Theorem \ref{thm:main}.
\qed

\bigskip \noindent
\textsc{Matan Karklinsky \\
Einstein Institute of Mathematics, The Hebrew University \\
Givat-Ram, Jerusalem 91904, Israel \\
Email: } \texttt{matan.karklinsky@mail.huji.ac.il}

\medskip
\bigskip \noindent
\textsc{Dan Romik \\
Einstein Institute of Mathematics, The Hebrew University \\
Givat-Ram, Jerusalem 91904, Israel \\
Email: } \texttt{romik@math.huji.ac.il}

\end{document}